# On some applications of complex numbers to polar equations and cycloidal curves

H. Azad, A. Laradji and M. T. Mustafa

Department of Mathematics & Statistics,

King Fahd University of Petroleum & Minerals, Dhahran 31261, Saudi Arabia

hassanaz@kfupm.edu.sa, alaradji@kfupm.edu.sa, tmustafa@kfupm.edu.sa

**Abstract**

The efficacy of using complex numbers for understanding geometric questions related to polar equations and general cycloids is demonstrated.

Complex numbers are nowadays taught in precalculus courses all over the world. Their integration in the teaching of classical topics related to Calculus is truly illuminating in every sense. We illustrate this by giving applications of complex numbers to the topics of polar equations and cycloidal curves. One needs the most basic ideas about complex numbers to use them effectively.

The topic of polar coordinates and polar curves – their sketching, points of intersection and related area problems – is often more challenging to teach to undergraduate students than its Cartesian counterpart (cf. [10, 11]).

One of the main reasons for this seems to be that not much emphasis is placed, in the traditional approach, on singular points of polar curves, whereas important geometric information is related to these points. Bringing such ideas to the fore is helpful in many problems related to polar curves, cf. [12].

One of the principal aims of this note is to demonstrate the great transparency – and consequent effectiveness - that results from using complex numbers in this context.

The other topic considered in this note is the theory of general cycloids. The theory of cycloids is an integral part of our cultural heritage. Yet this timeless and eminently practical mathematics is ignored or dealt with perfunctorily in the modern curriculum. The books of Brieskorn [3] and Moar [6] discuss this mathematics in great detail –its history, applications to technology and arts, as well as to the development of mathematics itself. In particular, Moar's book is at once a lament for the mathematics of a bygone era and a plea to restore these topics to the modern curriculum.

It is our opinion that these topics can be integrated in the modern curriculum and, in fact, are ideally suited for the use of technology in the classroom. Again, this can be done efficiently provided one is prepared to integrate complex numbers with calculus instruction.

Zwikker [8] has also used complex numbers systematically to study plane curves. The material on classical cycloidal curves is given in Chapter XIX of his book. In this paper, a more general situation is considered and the exposition is of the same level as for polar coordinates.

This note is written with our peers in mind .We have thus used concepts and terminology that we may not use in introductory courses and for which intuitive substitutes can always be found for classroom instruction. Also, some of the examples go beyond what would normally be taught in customary courses.

For reasons of clarity and exposition, we have only given details regarding concepts and definitions and left the verification of computations for the interested reader.

## 1 Polar equations

### 1.1 Intersection of polar curves

The advice in the text books we have seen, for example [2, 4, 5, 7] as well as in the Internet resources – for example [9]- is that one must draw the polar curves to find their points of intersection. This is misleading and goes against the spirit of our subject. In fact, finding points of intersection is no more or less difficult than solving equations. On the other hand, drawing curves defined by their polar equations is very difficult to do without technology and even if one sees a pair of curves using technology, their points of intersection still need to be computed to set up area integrals.

For the purposes of this note, a polar equation and its graph mean the following:

We take a level set of a function in the $(r, \theta)$ plane, say $f(r, \theta) = k$. Its image under the function $(r, \theta) \to (r\cos\theta, r\sin\theta)$ is its polar graph.

Let us consider the special case $r = f(\theta)$. Its polar graph is the set of all points $(x, y) = (f(\theta)\cos(\theta), f(\theta)\sin(\theta))$. Thinking of points as complex numbers and writing $(x, y) = z = x + iy$, the polar graph is the set of all points given by $z = f(\theta)e^{i\theta}$.

The main observation in computing points of intersections of polar graphs is this: if $r, s$ are non-zero real numbers and $\theta, \varphi$ are real then the equality $re^{i\theta} = se^{i\varphi}$ can only hold if $r = s$ and $\theta = \varphi \pmod{2\pi}$ or $r = -s$ and $\theta = \varphi + \pi \pmod{2\pi}$.

Thus, to find non-zero common points of the polar curves $r = f(\theta)$ and $r = g(\theta)$ one needs to solve the equation $f(\theta)e^{i\theta} = g(\varphi)e^{i\varphi}$, for non-zero numbers $f(\theta)$ and $g(\varphi)$.

This means that, in view of the observation above, one needs to find non-zero solutions of the equations $f(\theta) = g(\theta + 2n\pi), n\epsilon\mathbb{Z}$ and $f(\theta) = -g(\theta + \pi + 2n\pi), n\epsilon\mathbb{Z}$.

The origin has to be tested separately, as it has no well-defined polar coordinates-by checking if there is at least one $\theta$ and $\varphi$ for which $f(\theta) = 0$ and $g(\varphi) = 0$.

In particular if the functions $f$ and $g$ are $2\pi$ periodic, one needs only to solve for non-zero solutions of the equations $f(\theta) = g(\theta)$ and $f(\theta) = -g(\theta + \pi)$ and take care to test for the origin as a common point.

Before giving examples, a few words on symmetries and periodicity of polar curves are in order.

Calculus textbooks, such as the ones cited here, usually give tests that are used to check if a polar curve has one or more of the three standard symmetries: about the x-axis, the y-axis, and the origin. The disadvantage of these tests is that they may fail even when there is symmetry. For example, the curve $r = cos(\theta/2)$ has all three symmetries but changing $\theta$ to $\pi - \theta$ (for symmetry about the y-axis) does not yield an equivalent equation. We can obtain a more conclusive test once we note that a point with polar coordinates $(r, \theta)$ is on a polar curve $F(r, \theta) = 0$ precisely when there is an integer $n$ such that $F((-1)^n r, n\pi + \theta) = 0$.

Now, a rotation through an angle $\theta_0$ maps the point $re^{i\theta}$ to the point $re^{i(\theta+\theta_0)}$ while a reflection in the line $\theta = \theta_0$ maps the point $re^{i(\theta_0+\epsilon)}$ to the point $re^{i(\theta_0-\epsilon)}$ and therefore it maps the point $re^{i\theta}$ to the point $re^{i(2\theta_0-\theta)}$. Therefore using the equations that give equality of $re^{i\theta} = se^{i\varphi}$ for nonzero $r, s$, the reader can check that a nonzero point $f(\theta)e^{i\theta}$ of the polar curve defined by a function $r = f(\theta)$ is mapped to another point of the same curve by rotation through an angle $\theta_0$ if for some integer $n$– that could depend on $\theta$ – we have $f(\theta) = (-1)^n f(\theta + \theta_0 + n\pi)$, while a nonzero point $f(\theta)e^{i\theta}$ of the polar curve defined by a function $f$ is mapped to another point of the same curve by reflection through the line $\theta = \theta_0$ if, for some integer $n$– that could depend on $\theta$ – we have $f(\theta) = (-1)^n f(2\theta_0 - \theta + n\pi)$.

Symmetry about the origin and reflections along the coordinate axes now follow immediately from this.

We say that the polar graph defined by a real valued function $f$ is periodic with period $T$ if $T$ is the smallest positive real number such that $f(\theta)e^{i\theta} = f(\theta + T)e^{i(T+\theta)}$. Ignoring trivialities, this means $e^{iT}$ is a real number and therefore $T = N\pi$ for some positive integer $N$.

We illustrate our results with several examples.

**Example 1.** *Discuss the symmetries and periodicity of the curve with polar equation* $r = \cos\left(\frac{m\theta}{n}\right)$ *where* $\frac{m}{n}$ *is a positive rational number in its lowest terms.*

Since $\cos\left(\frac{m\theta}{n}\right) = (-1)^{2n} \cos\left(\frac{m(2n\pi - \theta)}{n}\right)$ it follows that the curve is symmetric with respect to the x-axis. If $m$ or $n$ is even, we have $\cos\left(\frac{m\theta}{n}\right) = (-1)^{n+1} \cos\left(\frac{m(n\pi - \theta)}{n}\right)$ and so the curve is also symmetric with respect to the y-axis and therefore the origin. If both $m$ and $n$ are odd then symmetry w.r.t. the y-axis would imply that for each $\theta$

$$\cos\left(\frac{m\theta}{n}\right) = (-1)^{h+1} \cos\left(\frac{m(h\pi - \theta)}{n}\right) \text{ for some integer } h.$$

Since, as the reader may check, the only values of $\theta$ for which this occurs are rational multiples of $\pi$, we obtain that the polar curve is not symmetric w.r.t. the y-axis (nor the origin) when $m$ and $n$ are odd.

Next, if $r = \cos\left(\frac{m\theta}{n}\right)$ has period $k\pi$ ($k$ positive integer) then the equation $\cos\left(\frac{m\theta}{n}\right) = (-1)^k \cos\left(\frac{m}{n}(\theta + k\pi)\right)$ implies that $\frac{km}{n}$ is an integer, and so $n$ divides $k$. A consideration of the parities of $m$ and $n$ then shows that the period is $2n\pi$ if $m$ or $n$ is even, and is $n\pi$ if $m$ and $n$ are odd.

In particular, to graph the equation $r = \cos\left(\frac{\theta}{N}\right)$ where $N$ is even, we can graph it for $0 \le \theta \le \frac{\pi N}{2}$ and use symmetry w.r.t. the coordinate axes to obtain the remaining part of the graph.

**Example 2.** *Find the number of points of intersection (apart from the origin) of the polar curves* $r = \cos(m\theta)$ *and* $r = \sin(n\theta)$ *where* $m$ *and* $n$ *are odd positive integers.*

Since the number of points of intersection is invariant if we interchange sine and cosine in the equations, we may assume that $m \ge n$. From Example 1 $r = \cos(m\theta)$ has period $\pi$ if $m$ is odd and $2\pi$ if $m$ is even. Suppose first that $m$ and $n$ are odd. A point of intersection $(r, \theta)$ must then satisfy $\cos(m\theta) = (-1)^h \sin(n(\theta + h\pi))$ for some integer $h$, i.e.

$$\cos(m\theta) = \sin(n\theta) = \cos\left(\frac{\pi}{2} - n\theta\right)$$

The required solutions satisfy $m\theta = \pm\left(\frac{\pi}{2} - n\theta\right) + 2k\pi$ where the integers $k$ are chosen so that $0 \le \theta < \pi$. Note that the minus sign cannot occur if $m = n$. We thus get two families of solutions satisfying: $0 \leq 2k + \frac{1}{2} < m+n$ i.e. $0 \leq k < \frac{m+n}{2}$ and, in case $m \neq n$, $0 \leq 2k - \frac{1}{2} < m - n$ i.e. $1 \leq k \leq \frac{m-n}{2}$. It is easy to see that points belonging to the same family are distinct. So, if a solution were repeated, there would be integers $h, k$ and $l$ such that $\frac{2k-\frac{1}{2}}{m-n}\pi + l\pi = \frac{2h+\frac{1}{2}}{m+n}\pi$. This is impossible modulo 2, and so the total number of points of intersection (apart from the origin) is $m$.

The next three examples are standard examples from textbooks.

**Example 3**. *Find all points of intersection of the curves with polar equations*

$r = 1$ and $r = cos\theta$

The intersection points are given by solutions of equations

$1 = cos\theta, \quad 1 = -\cos(\theta + \pi)$.

As $\cos(\theta + \pi) = -\cos(\theta)$, the points are given by just $cos\theta = 1$. There is only one solution modulo $2\pi$ and the intersection point is $z = 1$.

**Example 4.** *Find all points of intersection of the curves with polar equations*

$r = cos\theta, \quad r = 1 - cos\theta$.

The origin is a common point as the value 0 is in the range of both functions.
We can save some calculations by noticing that the second function is always non-negative. Thus the points other than the origin are given by solving

$1 - cos\theta = cos\theta$ and $1 - cos\theta = -\cos(\theta + \pi)$,

so, in fact we have to solve just one equation $1 - cos\theta = cos\theta$, which gives two points $z_1 = \frac{1}{2}e^{\frac{i\pi}{3}}, z_2 = \frac{1}{2}e^{-\frac{i\pi}{3}}$.

This, together with 0, gives 3 intersection points.

Working with the equations $cos\theta = 1 - cos\theta$ and $= -(1 - \cos(\theta + \pi))$ should, of course give the same points, as the reader can check.

**Example 5.** *Find all points of intersection of the curves with polar equations*

$r = 2cos\theta, r = 1 + cos\theta$

As in the previous example, one has just to solve the equation

$1 + cos\theta = 2cos\theta$

that is, $\cos\theta = 1$ to get the points other than the origin. This gives the two points of intersection 0 and 1.

**Example 6.** *Find all points of intersection of the curves with polar equations*
$$r = \sin N\theta, r = \cos N\theta,$$
*where N is a positive integer.*

Clearly 0 is a point of intersection. To find the other points of intersection, one needs to solve – as explained above – the equations
$$\sin N\theta = \cos N\theta \text{ and } \sin N\theta = -\cos N(\theta + \pi) = (-1)^{N+1}\cos(N\theta).$$
Thus for odd $N$ one needs to solve the equation $\tan N\theta = 1$ while for even $N$ one needs to solve the equations $\tan N\theta = 1$ as well as $\tan N\theta = -1$.

- **Case (i)** *N* is odd.

  The solutions of $\tan N\theta = 1$ are $\theta = \frac{\pi}{4N} + \frac{n\pi}{N}$, $n\epsilon\mathbb{Z}$.

  The nonzero points of intersection are thus
  $$p(n) = \sin\left(\frac{\pi}{4} + n\pi\right) e^{i\left(\frac{\pi}{4N} + \frac{n\pi}{N}\right)} = (-1)^n \left(\sin\frac{\pi}{4}\right) e^{\frac{i\pi}{4N}} e^{\frac{in\pi}{N}}, n\epsilon\mathbb{Z}.$$

  Thus to get all the points, we may assume that *n* takes values between 0 and 2*N*−1.

  The reader can check that if $p(l) = p(k)$ and $l, k$ have the same parity, then $l \equiv k \ (mod \ 2N)$; while if $l, k$ have opposite parities then necessarily $l \equiv k + N \ (mod \ 2N)$. Conversely, as $N$ is odd,

  if $l \equiv k + N \ (mod \ 2N)$ then $l, k$ have different parities and the corresponding points are the same.

  Thus, the distinct non-zero points of intersection of the curves are the $N$ points $p(0), ..., p(N-1)$

- **Case (ii)** *N* is even. Using the notation as in Case (i), the congruence $l \equiv k + N \ (mod \ 2N)$ forces $l, k$ to have the same parity. Thus we get 2*N* distinct points $p(0), ....., p(2N - 1)$ from the solutions of the equation $\tan N\theta = 1$.

  Similarly, using the solutions of $\tan N\theta = -1$, we get the distinct points
  $$q(n) = \sin\left(\frac{\pi}{4} + n\pi\right) e^{i\left(\frac{\pi}{4N} + \frac{n\pi}{N}\right)} = (-1)^n \left(\sin\frac{\pi}{4}\right) e^{\frac{i\pi}{4N}} e^{\frac{in\pi}{N}}, n = 0, ..., 2n - 1.$$

  It remains to check that the points $p(n), q(m), n, m = 0, ..., 2N - 1$ are distinct, which can be done again by working with congruences: the equality $p(n) = q(m)$ implies the impossible congruences $2n + 2N \equiv 1 + 2m (mod 4N)$ - in case $n, m$ have the same parity, and $2n \equiv 1 + 2m \ (mod \ 4N)$ in case $n$ is odd and $m$ is even.

Thus we get 4N distinct non-zero points of intersection if $N$ is even.

## 1.2 Computations relating to areas

In computing the area of intersection of the region bounded by polar curves, one has to be careful that in the equations defining the boundary curves $r = f(\theta)$, the function $f$ can be chosen so that it takes non-negative values. This can always be achieved by noticing that if a function $f$ is defined on an interval $I$ and all its values on $I$ are non-positive, then

$$f(\theta)e^{i\theta} = -f(\theta)e^{i(\theta+\pi)} = -f(\varphi - \pi)e^{i\varphi}, \varphi \epsilon I + \pi.$$

Thus we can replace the given function piecewise by non-negative valued functions to compute areas of intersection of regions.

The advantage of such a representation is that the non-zero points of intersection are much easier to calculate and the regions bounded by non-negative $2\pi$ periodic curves have simple descriptions. In particular, if the function $f$ is $2\pi$ periodic and the length of the interval $I$ is less than $2\pi$, then the only possible point of self-intersection is at the origin and the region bounded by $f$ is given by $0 \leq r \leq f(\theta), \theta \epsilon I$.

We illustrate this first by a very simple example. Here, it is certainly much easier to work with a sketch. But the example shows a certain point that is useful in more involved situations.

**Example 1:** *Find the area of intersection of the region bounded by the curves with polar equations $r = sin\theta, r = cos\theta$.*

As $r = sin\theta$ is negative in the interval $[\pi, 2\pi]$, the points represented by this equation in this interval are $- sin(\varphi - \pi) e^{i\varphi} = sin\varphi e^{i\varphi}, \varphi \epsilon [2\pi, 3\pi]$, so by periodicity of the sine function, the curve is traced twice and we need to only consider the curve in the interval $[0, \pi]$.

Similarly, for the points traced by $r = cos\theta$ we can restrict $\theta$ to the interval $[\frac{-\pi}{2}, \frac{\pi}{2}]$.

The regions bounded by these curves are defined by

$$0 \leq r \leq sin\theta, 0 \leq \theta \leq \pi, 0 \leq r \leq cos\theta, \frac{-\pi}{2} \leq \theta \leq \frac{\pi}{2},$$

so the intersection is given by $0 \leq r \leq \min(\sin\theta, \cos\theta), 0 \leq \theta \leq \frac{\pi}{2}$.

We need to solve the equation $\sin\theta = \cos\theta$ in this interval and that gives the area as a sum of two integrals.

**Example 2.** *Find the area of the intersection of the region bounded by the curves with polar equations $r = \sin N\theta$, $r = \cos N\theta$, where N is a positive integer*

As the functions are $2\pi$ periodic and 0 is in their range, the only point of self-intersection is at the origin.

Also, $\sin N\theta = 0$ for all multiples of $\frac{\pi}{N}$ and $\cos N\theta = 0$ for odd multiples of $\frac{\pi}{2}$.

First consider the polar graph of $r = \sin N\theta$. We get $2N$ tangent rays at the origin for this graph and the sign changes alternately.

Let $I$ be an interval on which $\sin N\theta$ is positive except at the end points, where it has a zero, say $I = [0, \frac{\pi}{N}]$. Rotation by $\frac{\pi}{N}$ maps the points $(\sin N\theta)e^{i\theta}, \theta \in I$ to the points

$$(\sin N\theta)e^{i\left(\theta + \frac{\pi}{N}\right)} = (-1)^N (\sin N\varphi)e^{i\varphi}, \varphi \in I + \frac{\pi}{N}.$$

Thus, for even $N$ the graph is preserved, whereas for odd $N$ the graph is mapped to the graph of $r = -\sin N\theta$.

Similar considerations apply to the graph of $r = \cos N\theta$.

Since the sine and cosine functions are positive in the first quadrant, then, as explained just before Example 1, to compute the required area, we need only to compute the area common to the region $0 \leq r \leq \sin N\theta, 0 \leq \theta \leq \frac{\pi}{N}$ and the region $0 \leq r \leq \cos N\theta, \frac{-\pi}{2N} \leq \theta \leq \frac{\pi}{2N}$ and multiply it by $N$ or $2N$, depending whether $N$ is odd or even.

The common area is described in this sector by

$$0 \leq r \leq \min(\sin N\theta, \cos N\theta), 0 \leq \theta \leq \frac{\pi}{2N}$$

and thus by

$$\left\{(r,\theta): 0 \leq r \leq \cos N\theta, 0 \leq \theta \leq \frac{\pi}{4N}\right\} \cup \left\{(r,\theta): 0 \leq r \leq \sin N\theta, \frac{\pi}{4N} \leq \theta \leq \frac{\pi}{2N}\right\}.$$

The next example in this section is suggested by a question in the Physics Forum.

**Example 3** *Limaçons*

*Find the area of the region inside both the large loop of the curve $r = 1 - \lambda \sin\theta$ and the small loop of the curve $r = 1 + \lambda \cos\theta$, where $\lambda > 1$.*

As in the previous examples, all essential geometric information is contained in the points of self-intersection.

We need to describe these limaçons by piecewise positive functions.

First consider the limaçon with polar equation

$r = 1 - \lambda \sin\theta$. Now $r = 0$ at $\theta_0 = \sin^{-1}(\frac{1}{\lambda})$ and at $\pi - \theta_0$.

The maximum value of $r$ is $\lambda + 1$ and the minimum is $1 - \lambda$.

The large loop is therefore given by

$r = 1 - \lambda \sin\theta, \quad \pi - \theta_0 \leq \theta \leq \theta_0 + 2\pi$ .

and the small loop by

$r = -(1 + \lambda \sin\theta), \quad \theta_0 + \pi \leq \theta \leq \pi - \theta_0 + \pi$

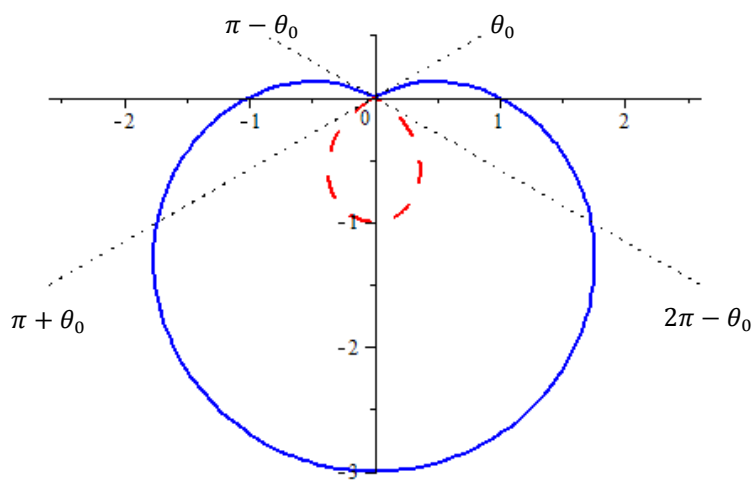

Figure 1

Similarly, for the polar equation $= 1 + \lambda \cos\theta$ , with $\lambda>1$, let $\varphi_0 = \cos^{-1}\left(-\frac{1}{\lambda}\right)$. The equations for the large loop and the small loop are respectively

$r = 1 + \lambda \cos\theta , \quad -\varphi_0 \leq \theta \leq \varphi_0$

and

$r = (\lambda\cos\theta - 1), \varphi_0 + \pi \leq \theta \leq -\varphi_0 + \pi + 2\pi.$

Moreover $\varphi_0 = \frac{\pi}{2} + \theta_0$ .

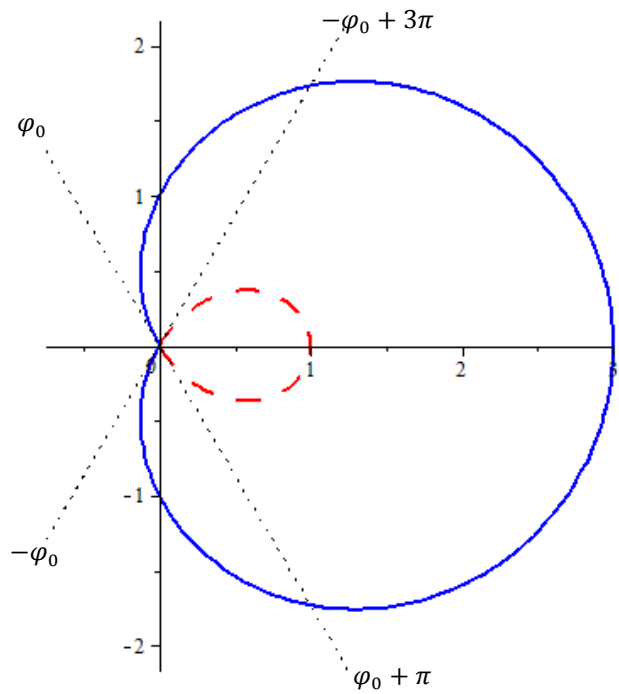

Figure 2

Thus the large loop of the first limaçon is given by
$$r = 1 - \lambda \sin\theta, \quad \pi - \theta_0 \le \theta \le \theta_0 + 2\pi$$
and the small loop of the second limaçon is given by
$$r = (\lambda \cos\theta - 1), \quad \theta_0 + \frac{\pi}{2} + \pi \le \theta \le -(\theta_0 + \frac{\pi}{2}) + 3\pi.$$

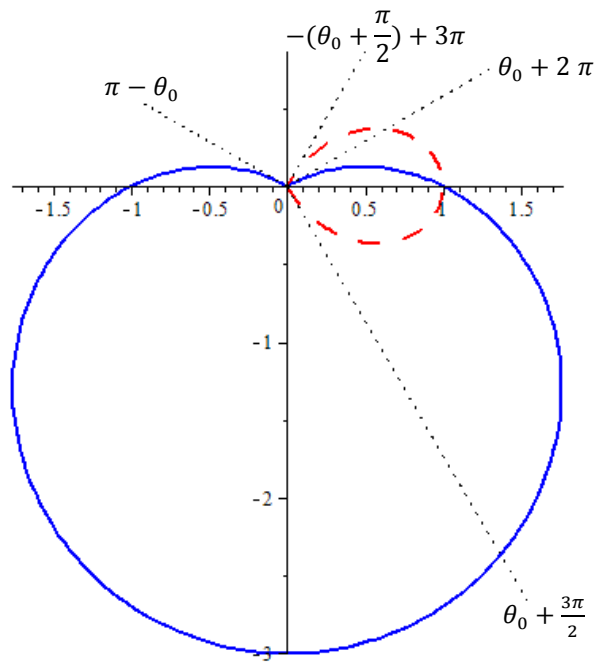

Figure 3

As these loops are now defined by positive functions, the intersection can be computed easily.

Clearly $\pi - \theta_0 < \theta_0 + \frac{\pi}{2} + \pi < \theta_0 + 2\pi$ and the right end-point of the interval for the small loop is in this same interval if and only if $\frac{\pi}{4} \leq \theta_0$. In this case, $\frac{1}{\sqrt{2}} \leq \frac{1}{\lambda}$, and the entire small loop is contained in the large loop.

So suppose $\frac{\pi}{4} > \theta_0$. Then, the intersection is in the sector given by

$$\theta_0 + \frac{3\pi}{2} < \theta < \theta_0 + 2\pi.$$

The point of intersection of the curves is given by $\lambda \cos\theta - 1 = 1 - \lambda \sin\theta$. Thus it is given by $\lambda(\cos\theta + \sin\theta) = 2$, that is by $\sin\left(\theta + \frac{\pi}{4}\right) = \frac{\sqrt{2}}{\lambda}$. This has two solutions, say $\theta_1, \theta_2$, where $\theta_1 + \frac{\pi}{4} = \sin^{-1}(\frac{\sqrt{2}}{\lambda})$ and $\theta_2 + \frac{\pi}{4} = \pi - (\theta_1 + \frac{\pi}{4})$.

Therefore $\theta_2 = \frac{\pi}{2} - \theta_1$.

There is only one non-zero point of intersection of the curves: the solution $\theta_1$ is between $-\frac{\pi}{4}, \frac{\pi}{4}$ and $\theta_2$ is outside the sector of intersection. Calculations show that the solution $\theta_1$ is given by $\sin^{-1}(\frac{1}{\lambda} - \sqrt{\frac{1}{2} - \frac{1}{\lambda^2}})$.

**Example 4** *2Nπ-periodic functions*

The polar graph of a $2\pi N$ − periodic function is the same as the graphs − each over the interval $[0, 2\pi]$ − of the functions

$$r = f(\theta), r = f(\theta + 2\pi), \ldots, r = f(\theta + (N-1)2\pi).$$

In problems involving areas, one has to describe these graphs by positive functions. If $f$ and $g$ are functions defined over an interval $I$ and $f(\theta) = g(\theta)$ for some $= \theta_1$, then the sign of the difference $f(\theta) - g(\theta)$ near $\theta_1$ can be determined from the finite Taylor series of $f(\theta) - g(\theta)$ of appropriate order.

For example, the polar graph of $r = \cos\left(\frac{\theta}{2}\right)$ is the same as the graphs of $r = \cos\left(\frac{\theta}{2}\right)$ and $r = -\cos\left(\frac{\theta}{2}\right)$ − both over the interval $[0, 2\pi]$. The second graph is obtained from the first by rotating it through 180-degrees.

The description of the graph of $r = \cos\left(\frac{\theta}{2}\right)$ over the interval $[0, 2\pi]$ as graphs of positive functions is given by the graphs of $r = \cos\left(\frac{\theta}{2}\right)$ and $r = \sin\left(\frac{\theta}{2}\right)$ over the interval $[0, \pi]$. The reason is that if $\theta \epsilon [\pi, 2\pi]$ then $\cos\left(\frac{\theta}{2}\right)$ is negative and therefore, writing $\theta = \varphi + \pi$, we have

$$\cos\left(\frac{\theta}{2}\right)e^{i\theta} = \sin\left(\frac{\varphi}{2}\right)e^{i\varphi}, \varphi \epsilon [0, \pi].$$

These curves intersect at

$$\theta = \varphi = \pi/2$$

and this gives the area of the intersection of the regions bounded by these two positive valued polar curves as a sum of two integrals.

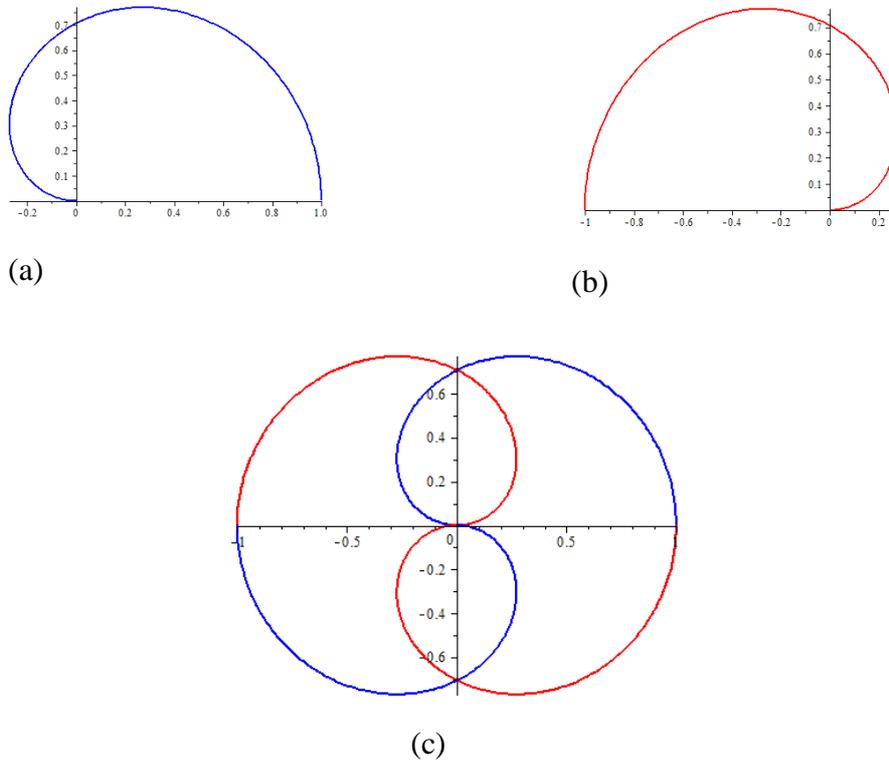

(a)

(b)

(c)

Figure 4. (a) Graph of $r = \cos\left(\frac{\theta}{2}\right)$ over the interval $[0, \pi]$. (b). Graph of $r = \sin\left(\frac{\theta}{2}\right)$ over the interval $[0, \pi]$. (c) Graph of $r = \cos\left(\frac{\theta}{2}\right)$ over the interval $[0, 4\pi]$ is obtained by combining Figures 1(a) and (b) and rotating the resulting graph by 180 degrees.

Similar considerations apply to periodic functions whose periods are odd multiples of π.

## 2 Cycloidal curves

These are the curves which are obtained by tracing the initial point of contact of a circle as it rolls without slipping on a line or another circle.

The books [3, 6, 8] discuss all aspects of these curves in great detail and contain beautiful applications.

Our aim here is only to derive the parametric equations of the curve traced by the initial point of contact of a circle as it rolls without slipping on another curve – either inside or outside the region bounded by the curve- whenever the concept of inside or outside a region makes sense -and then derive the equations of the cycloid, epicycloid and hypocycloid as special cases.

A more general situation can be understood by the same ideas, as explained in the remark immediately before example 1 of this section.

Instead of getting involved with what it means to be inside or outside a region, it is more efficient to deal with parameterized curves and derive the equations as a circle rolls on the curve without slipping. This can be done in remarkably straightforward manner if one is prepared to use complex numbers. The reason is that it is very easy to deal with planar rotations using complex numbers.

The proof is a translation of the following figures in a more general context:

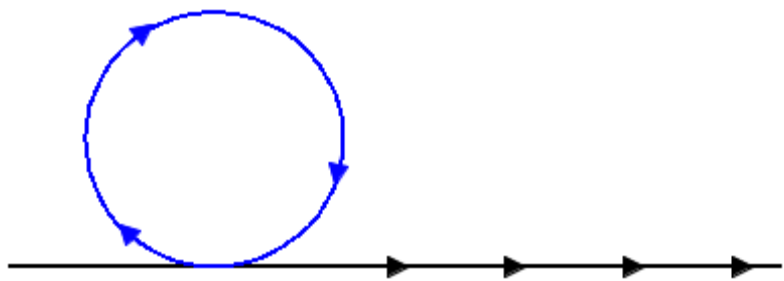

Figure 5

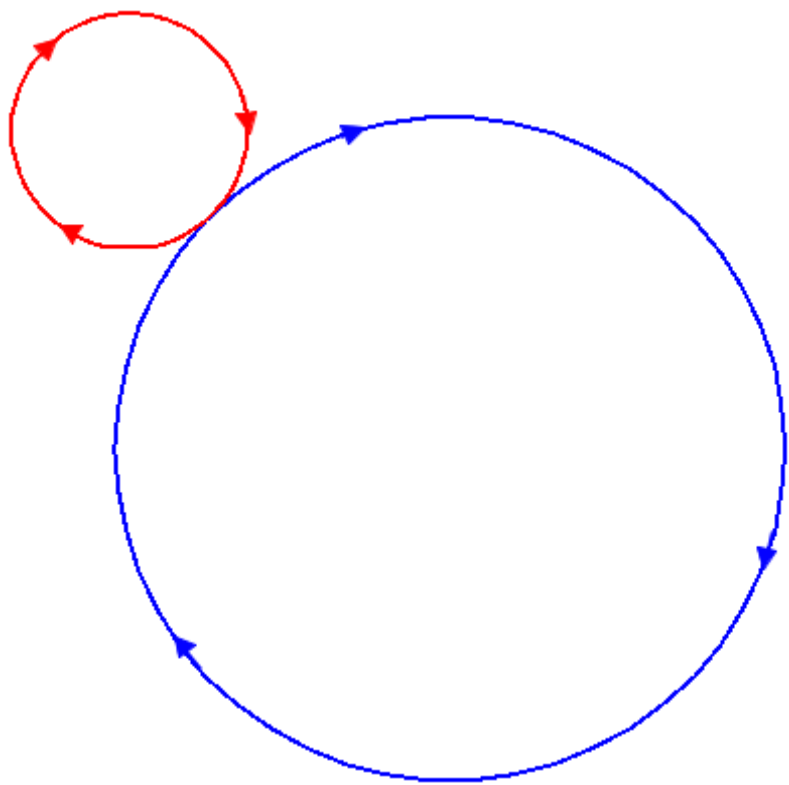

Figure 6

Let $\alpha: I \to \mathbb{R}^2$ be a regular parameterized curve: this means that the tangent vector $\alpha'(t)$ at the point $\alpha(t)$ is never zero. There are then two normal vectors at the point for the parametric value $t$ - namely those obtained by rotating $\alpha'(t)$ clockwise and anticlockwise through 90 degrees. Since rotation by 90 degrees anticlockwise is given by multiplication by the complex number $i$, this gives us the principal normal field $n(t)$ along the curve, namely $n(t) = i\alpha'(t)$.

We have a circle of radius $r$ that starts to roll tangentially above the curve $\alpha$ from left to right, without slipping. P is the initial point of contact. We want to describe the trajectory of the point P in this motion; as shown in Figure 7:

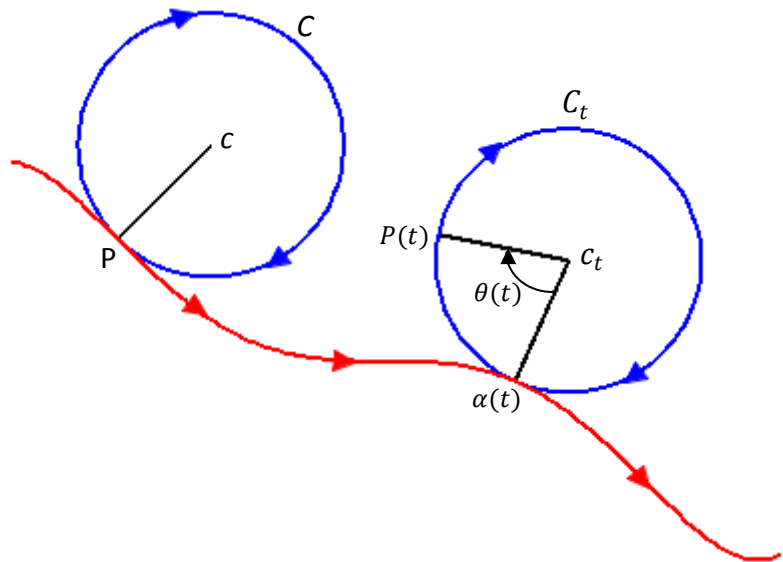

Figure 7

But what does left and right, above and below etc. mean for general curves? This can be easily extracted from Figure 7.

To stay close to intuition, let the parameter $t$ represent time.

To say that a circle $C$ of radius $r$ rolls on the curve $\alpha$ means we have a family of circles $C_t$ of radius $r$ that are all tangential to the curve $\alpha$, such that at time $t$, the center $c_t$ of the circle $C_t$ is on the normal line at the point $\alpha(t)$.

There are several possibilities. We discuss only two that represent physical motion.

The direction of the vector $\overrightarrow{\alpha(t)c_t}$ is either that of the principal normal $n(t) = i\alpha'(t)$ or that of the vector $-n(t)$.

In the first case, we assume that the circle is rotating clockwise and in the second case that it is rotating anticlockwise. In both cases, we assume that the motion is without slippage.

More precisely, this means, in the first case that

(i) The initial point of contact $P = \alpha(t_0)$ has moved to the point $P(t)$ on $C_t$ obtained by rotating clockwise the vector $\overrightarrow{c_t\alpha(t)}$ at $c_t$ through an angle $\theta(t)$, which is determined by the following requirement:

(ii) the distance traveled on the circle by the initial point of contact $P$ to the position $P(t)$ on the circle is the same as the distance along the curve from the point $P = \alpha(t_0)$ to the point of contact $\alpha(t)$ at time t.

This means that $r\theta(t)$ equals the distance traversed along the curve $\alpha$ from $\alpha(t_0)$ to $\alpha(t)$.

Translating this into equations we have:

$$c_t = \alpha(t) + i \frac{\alpha'(t)}{|\alpha'(t)|} r$$

$$\overrightarrow{c_t P(t)} = \overrightarrow{c_t\alpha(t)}\, e^{-i\theta(t)}$$

$$P(t) = c_t + \overrightarrow{c_t P(t)}$$

$$= \alpha(t) + i \frac{\alpha'(t)}{|\alpha'(t)|} r - ir\frac{\alpha'(t)}{|\alpha'(t)|} e^{-i\theta(t)} \qquad (*)$$

$$r\theta(t) = \int_{t_0}^{t} |\alpha'(u)|\, du$$

Similarly, in the second case, we assume that the initial point of contact $P = \alpha(t_0)$ has moved to the point $P(t)$ on $C_t$ obtained by rotating anti-clockwise the vector $\overrightarrow{c_t\alpha(t)}$ at $c_t$ through an angle $\theta(t)$, which is determined by the equality of distance traveled on the circle by the initial point of contact and the distance along the curve from the initial point of contact to the point $\alpha(t)$ of contact at time t.

Changing every plus sign into minus and vice-versa in the equation (*), we obtain the equations for the position of the initial point of contact when a circle rolls as described in the second case.

We record this information in the following equation:

$$P(t) = c_t + \overrightarrow{c_t P(t)}$$

$$= \alpha(t) - i \frac{\alpha'(t)}{|\alpha'(t)|} r + i \frac{\alpha'(t)}{|\alpha'(t)|} e^{i\theta(t)} r \qquad (**)$$

$$r\theta(t) = \int_{t_0}^{t} |\alpha'(u)|\, du$$

The equations (*) and (**) are the main equations of this section. Specializing to straight lines and circles we get the equations of the cycloid, epicycloid and the hypocycloid.

*Remark:* The same ideas can be used to describe trajectory of the points *Q(t)* determined by the following conditions.

- *Q(t)* is on the line joining *P(t)* and $c_t$.
- The vector $\overrightarrow{P(t)Q(t)} = k(t)\overrightarrow{c_tP(t)}$ where $k(t)$ is a function of *t*. Usually one takes $k(t)$ to be a constant.

**Example 1.** *The parametric equations of a cycloid*

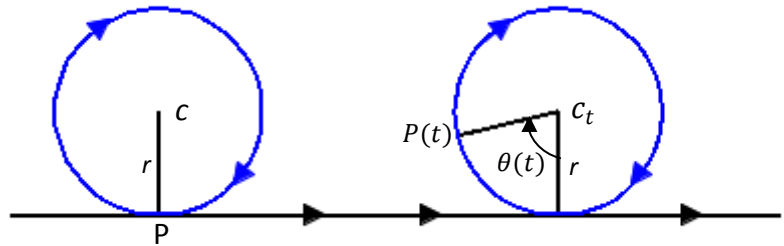

Figure 8

This is the curve traced out by a point on a circle of radius $r$ as it rolls on a straight line without slipping. We take the line to be the *x*-axis and think of it as an oriented curve by setting $\alpha(t)=t$. The principal normal at $t$ is given by $i$, so we should use the equation (*) above. The parametric equations of the cycloid are therefore given by $P(t) = \alpha(t) + i\frac{\alpha'(t)}{|\alpha'(t)|}r - i\frac{\alpha'(t)}{|\alpha'(t)|}e^{-i\theta(t)}r$ (*)

where

$$r\theta(t) = \int_0^t |\alpha'(u)| du$$

so they are given by

$$P(t) = t + ir - ire^{-i\frac{t}{r}}.$$

Thus

$$x(t) = \left(t - r\sin\left(\frac{t}{r}\right)\right), \quad y(t) = r - r\cos\left(\frac{t}{r}\right)$$

**Example 2.** *The parametric equations of an epicycloid.*

This is the curve generated by the trajectory of the initial point of contact of a circle of radius r rolling on top of a circle of radius $R$, with $R > r$.

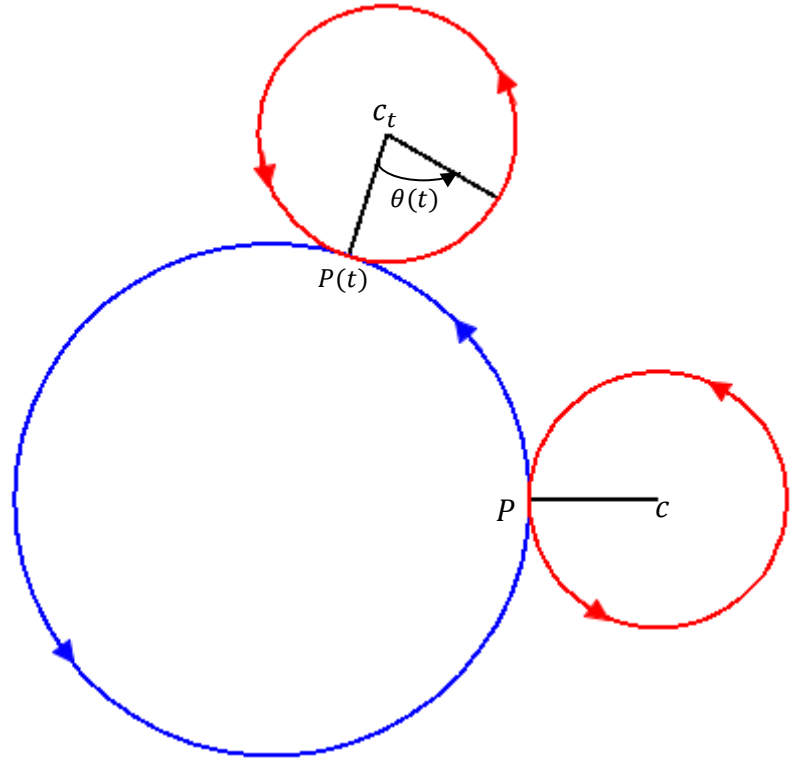

Figure 9

For the usual parameterization of the circle given by $\alpha(t) = Re^{it}$, it is traced counter-clockwise and so the principal normal points inwards. Therefore the motion of the small circle is modeled by the second case discussed above. The trajectory of the initial point of contact of the small circle as it rolls on top of the large circle is given by the equations (**):

$$P(t) = \alpha(t) - i \frac{\alpha'(t)}{|\alpha'(t)|} r + i \frac{\alpha'(t)}{|\alpha'(t)|} e^{i\theta(t)} r,$$

where

$$r\theta(t) = \int_0^t |\alpha'(u)| du.$$

This works out to

$$P(t) = (R + r)e^{it} - r\, e^{it(1+\frac{R}{r})}$$

which gives the standard parameterization

$$x(t) = (R + r)cost - rcost\left(1 + \frac{R}{r}\right),$$

$$y(t) = (R + r)\sin t - r\sin(1 + \tfrac{R}{r})$$

**Example 3.** *Parametric equations of a hypocycloid.*

This is the curve generated by the trajectory of the initial point of contact of a small circle of radius r rolling inside a circle of radius $R$, with $R > r$.

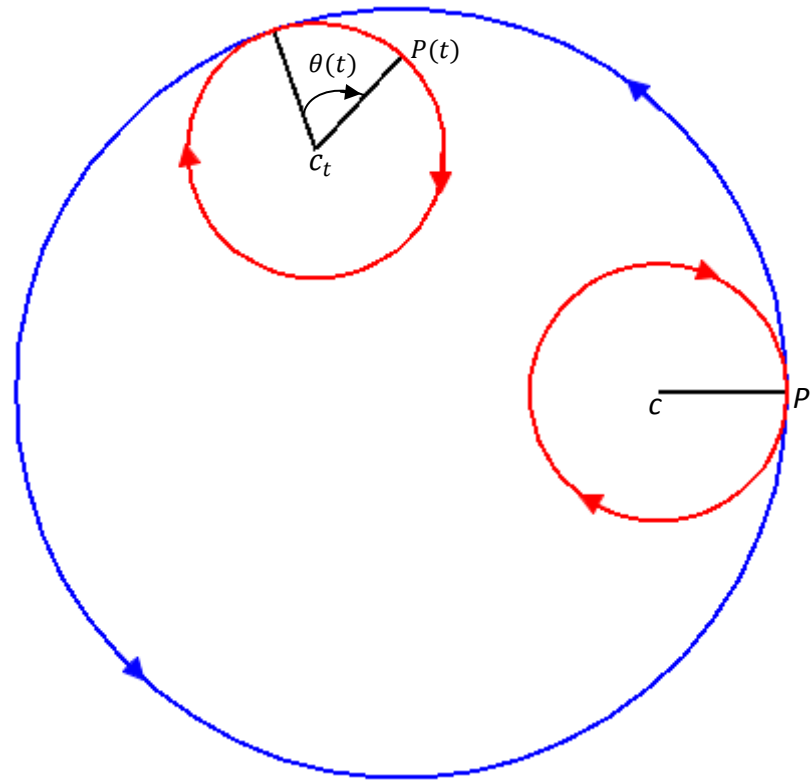

Figure 10

As the principal normal in the usual parameterization of the circle given by

$$\alpha(t) = Re^{it}$$

points inwards, the motion of the small circle is modeled by the first situation discussed above.

The equations of the hypocycloid are therefore given by

$$c_t = \alpha(t) + i \frac{\alpha'(t)}{|\alpha'(t)|} r$$

$$P(t) = \alpha(t) + i \frac{\alpha'(t)}{|\alpha'(t)|} r - ir \frac{\alpha'(t)}{|\alpha'(t)|} e^{-i\theta(t)}$$

where

$$r\theta(t) = \int_0^t |\alpha'(u)| du$$

This works out to

$$P(t) = (R-r)e^{it} + r\, e^{it(1-\frac{R}{r})}$$

which gives the parameterization

$$x(t) = (R-r)\cos t + r\cos t\left(1 - \frac{R}{r}\right),$$
$$y(t) = (R-r)\sin t + r\sin t\left(1 - \frac{R}{r}\right)$$

These simply generated curves have remarkable properties. The reader is referred to the books [3, 6, 8] for authoritative accounts of the subject.

Next we give examples of some curves generated by motions of circles on limaçons.

**Example 4.** *Motion of circle on limaçons.*

The parametric equations of the curves generated by the trajectory of the initial point of contact of a circle of radius $r$ rolling on top of a limaçon or inside a limaçon can be directly obtained using equations (**) or (*) respectively.

Figures 11 and 12 give illustrations of some curves generated by motion of circle on limaçons.

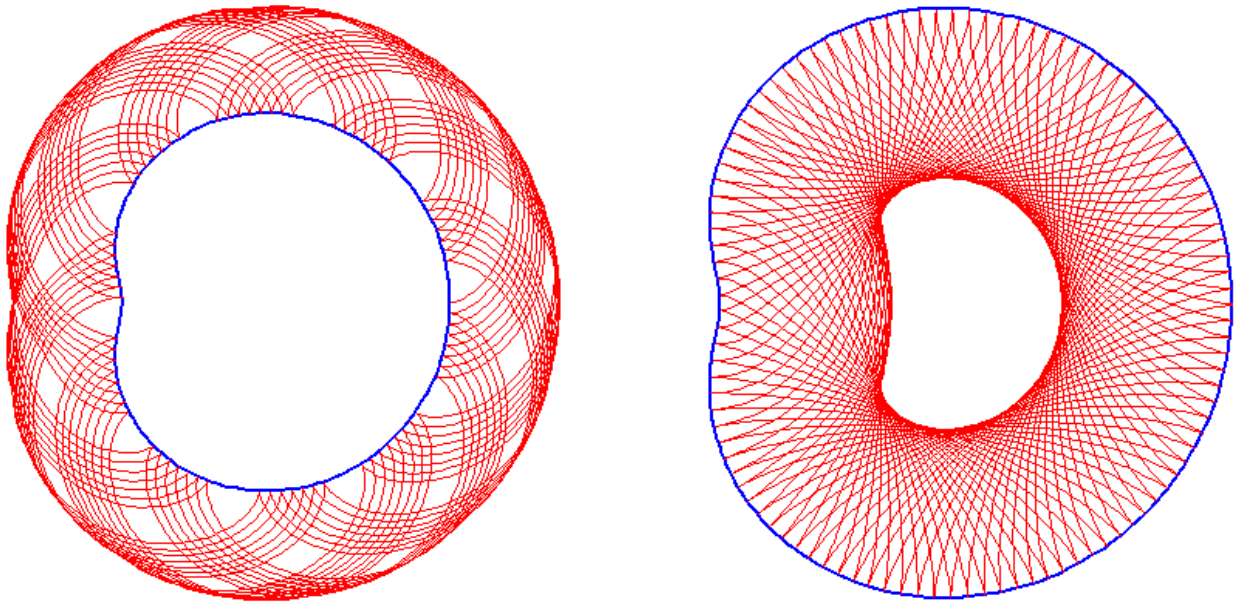

Figure 11

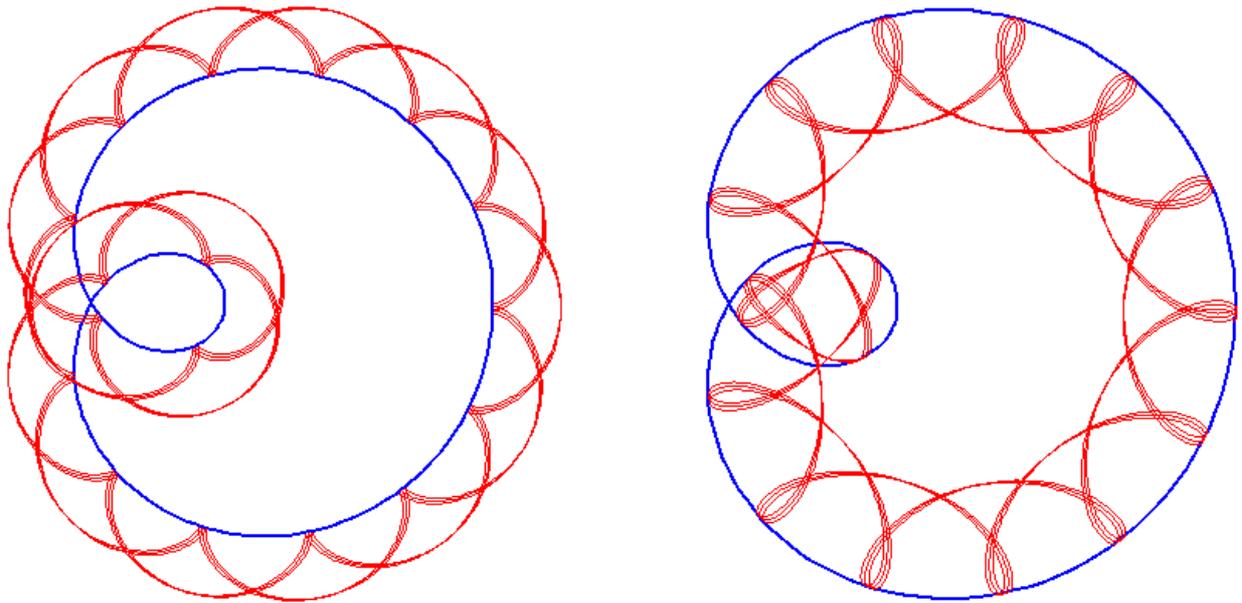

Figure 12

Reversing the configurations in the above analysis, that is, replacing $\theta(t)$ by $-\theta(t)$ in equations (*) and (**) can also be utilized to produce examples of curves that may not be easy to obtain through the motion with standard configuration. Illustrations of curves generated through the motion of circle, with reverse configuration, on circle and ellipse are provided in the next example. The equations for the motion of circle, with standard configuration, on an ellipse were also derived in [1] by a different approach.

**Example 5.**

Figure 13 shows examples of curves generated through the motion of a circle, with reverse configuration, rolling on top of a circle: here one gets the classical curves in Brieskorn [3, p.22].

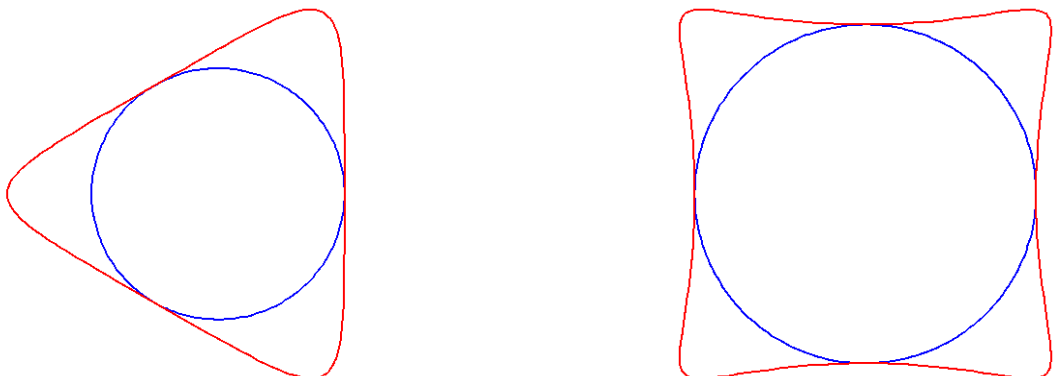

Figure 13

An example of a curve generated through the motion of circle, with reverse configuration, rolling inside an ellipse is given in Figure 14.

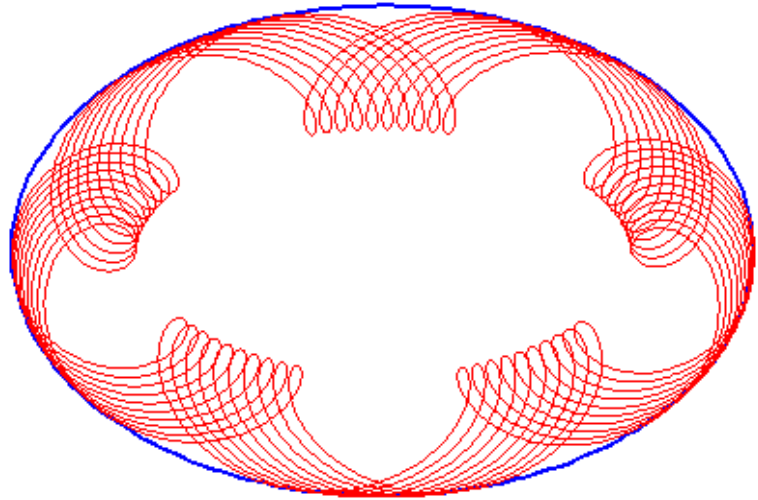

Figure 14

We hope that we have succeeded in showing that integrating the most basic facts about complex numbers in the teaching of classical topics is illuminating as well as very efficient.